# PARAMETER TUNING IN POINTWISE ADAPTATION USING A PROPAGATION APPROACH

By Vladimir Spokoiny and Céline Vial

*Weierstrass-Institute and Humboldt University Berlin, Modal'X and CREST-ENSAI, Modal'X, IRMAR*

This paper discusses the problem of adaptive estimation of a univariate object like the value of a regression function at a given point or a linear functional in a linear inverse problem. We consider an adaptive procedure originated from Lepski [*Theory Probab. Appl.* **35** (1990) 454–466.] that selects in a data-driven way one estimate out of a given class of estimates ordered by their variability. A serious problem with using this and similar procedures is the choice of some tuning parameters like thresholds. Numerical results show that the theoretically recommended proposals appear to be too conservative and lead to a strong oversmoothing effect. A careful choice of the parameters of the procedure is extremely important for getting the reasonable quality of estimation. The main contribution of this paper is the new approach for choosing the parameters of the procedure by providing the prescribed behavior of the resulting estimate in the simple parametric situation. We establish a non-asymptotical "oracle" bound, which shows that the estimation risk is, up to a logarithmic multiplier, equal to the risk of the "oracle" estimate that is optimally selected from the given family. A numerical study demonstrates a good performance of the resulting procedure in a number of simulated examples.

**1. Introduction.** This paper discusses the problem of selecting one estimate from a given family of estimates $\{\widetilde{\theta}_k, k = 1, \ldots, K\}$ of a univariate object $\theta$. We suppose that every estimate can be represented as

(1.1) $$\widetilde{\theta}_k = \theta_k + \xi_k, \qquad k = 1, \ldots, K,$$

where $\theta_k$ is the expectation of $\widetilde{\theta}_k$: $\mathbf{E}\widetilde{\theta}_k = \theta_k$ and $\xi_1, \ldots, \xi_K$ are zero mean random errors. In what follows we assume that $(\xi_1, \ldots, \xi_K)$ is a Gaussian vector with a known covariance matrix $B$. This problem is illustrated by two major examples: estimating a regression function at a given point and estimating a linear functional in a linear inverse problem. In the case of a Gaussian regression model $Y_i = f(X_i) + \varepsilon_i$, the target of estimation is the value of the unknown regression function $f(x)$ at a certain point $x$. The set $\{\widetilde{\theta}_k\}$ can be obtained as kernel or local polynomial estimates with different bandwidths. In the case of a linear inverse problem, the target is usually the value of a linear functional and the family of estimates is obtained by using different values of the regularization parameter for the regularized inversion. Note that the representation (1.1) can be regarded as a reasonable approximation for many other statistical models and problems like regression with non-Gaussian errors or density estimation.

The problem of adaptive estimation can be formulated as the best possible choice of one estimate out of this family on the basis of the available information. This problem can be viewed as the problem of model selection, see, for example, Birgé and Massart (1993, 1998), Birgé (2006), Juditsky, Rigollet and Tsybakov (2008) and references therein. However, there is an essential difference between the (global) model selection problem and the problem of pointwise estimation considered in this paper. In the problem of global model selection one tries to recover the whole underlying model, that is, the target is the model itself. Here we consider the problem of recovering a one-dimensional characteristic of the whole model like the value of the function at a certain point. This makes these two problems quite different. In particular, for the problem of pointwise adaptation some additional assumptions on the considered family of estimates are required. Typically one assumes that the given family of "weak" estimates $\widetilde{\theta}_k$ is ordered in the sense that the variance $v_k$ of $\widetilde{\theta}_k$ decreases with $k$. Another intrinsic assumption on the considered set-up is that the squared bias $b_k^2 \stackrel{\text{def}}{=} (\theta_k - \theta)^2$ is small for the $k = 1$ and it may increases with $k$. The most popular example is given by kernel estimates with different bandwidths so that the starting bandwidth $h_1$ is small leading to the small bias but a huge variance of estimation. As the bandwidth grows the variance decreases but the bias may increase dramatically. The aim is to construct from the data one estimate that performs in the best possible way and particularly minimizes the corresponding estimation risk.

The first adaptive procedure of this sort was suggested in Lepski (1990) and extended in Lepski (1992) to much more general set-up. The idea is to select the largest index $k$ such that the estimates $\widetilde{\theta}_1, \ldots, \widetilde{\theta}_k$ do not differ significantly with each other. Two estimates $\widetilde{\theta}_l$ and $\widetilde{\theta}_k$ for $l < k$ differ significantly if the standardized difference $T_{lk} \stackrel{\text{def}}{=} v_l^{-1}(\widetilde{\theta}_l - \widetilde{\theta}_k)^2/2$ exceeds the prescribed threshold $\mathfrak{z}$, which can be dependent of $l$, $\mathfrak{z} = \mathfrak{z}_l$. Lepski (1990) stated the rate optimality of this procedure over Hölder smoothness classes, and Lepski, Mammen and Spokoiny (1997) showed its spatial adaptivity in the sense of rate optimality over Besov functional classes and established some oracle risk bound. Lepski and Spokoiny (1997) proved sharp optimality of a slightly modified procedure in the asymptotic minimax sense. However, all the mentioned results have been established under some conditions on the thresholds $\mathfrak{z}_k$, which basically means that the thresholds have to be sufficiently large, and they tell nothing if this condition is not fulfilled. At the same time, numerical results in simulated and practical data examples show that applying a large threshold typically leads to a conservative procedure and oversmoothing effects. In this sense, one can say there is some critical gap between the theory and practical applications.

Our paper presents a novel method for selecting the tuning parameters of the method based on the so called "propagation" condition, which postulates the desirable performance of the method in the simple parametric situation. The idea is similar to the problem of hypothesis testing for which the critical value of a test is selected by bounding the first-kind error probability under the null hypothesis. The theoretical study is done for the adaptive estimate with the selected tuning parameters. The main result claims the desired oracle risk bound for this defined procedure. The proposed approach seems to be quite general and it can be directly applied to many other procedures including local model selection, stagewise aggregation and local change-point analysis, which are studied in details in Spokoiny (2009) in a much more general set-up.

Golubev (2004) proposed another "risk envelope" approach to select the threshold for a special sequence space model and a particular linear functional. We consider this example in Section 1.3. The common point between Golubev (2004) and our proposal is the selection of the parameters of the method by a Monte Carlo simulation from the model with zero response. However, the procedure, motivation and theoretical analysis of our study is quite different from the one in Golubev (2004).

Theoretical properties of the proposed method are presented in Section 3. The main result states the "oracle" property of the proposed estimate: the risk of the adaptive estimate is within a log-multiple as small as the risk of the "oracle" estimate for the given model. The results are established in the precise nonasymptotic way in a rather general form. Our simulation study in Section 4 confirms a nice finite sample performance of the procedure for a rather big class of different models and problems.

Below in this section we present three major examples for which the proposed procedure can be applied. We start with the problem of pointwise bandwidth selection in kernel estimation, then we discuss the problem of estimating a linear functional in a linear inverse problem and then specify it to one particular functional in the sequence space model.

1.1. *Bandwidth choice in kernel estimation.* Consider a regression model $Y_i = f(X_i) + \varepsilon_i$ where $\varepsilon_i$ are i.i.d. Gaussian errors with zero mean and known variance $\sigma^2$ and with a deterministic design $X_1, \ldots, X_n$ in $\mathbb{R}^d$. The considered problem is to estimate the value of the unknown regression function $f(x)$ at a given point $x$. Let a sequence of localizing scheme $W^{(k)} = \{w_i^{(k)}\}$ have been fixed for $k = 1, \ldots, K$. In the case of kernel weights, this sequence is built just by using different values of the bandwidth $h$ from the smallest bandwidth $h_1$ to the largest value $h_K$ in the form $w_i^{(k)} = \psi(|X_i - x|/h_k)$ with a kernel function $\psi(\cdot)$. Every localizing scheme yields the corresponding estimate

$$\widetilde{\theta}_k = N_k^{-1} \sum_{i=1}^n w_i^{(k)} Y_i, \qquad N_k = \sum_{i=1}^n w_i^{(k)}.$$

By simple algebra

$$\theta_k \stackrel{\text{def}}{=} \mathbf{E}\widetilde{\theta}_k = N_k^{-1} \sum_{i=1}^n w_i^{(k)} f(X_i), \qquad \xi_k \stackrel{\text{def}}{=} \widetilde{\theta}_k - \mathbf{E}\widetilde{\theta}_k = N_k^{-1} \sum_{i=1}^n w_i^{(k)} \varepsilon_i.$$

Moreover,

$$\mathbf{E}\xi_l \xi_k = \sigma^2 N_k^{-1} N_l^{-1} \sum_{i=1}^n w_i^{(k)} w_i^{(l)}.$$

The above ordering condition can be written for the case of the kernel weights in the form $N_l < N_k$ for $l < k$. Below we will assume even a stronger condition that the values $N_k$ grow exponentially.

1.2. *Estimation of a linear functional in a linear inverse problem.* Consider a general set-up of a linear inverse problem when the observed data $Y$ from a Hilbert space $\mathcal{H}_Y$ are modeled by a linear operator equation

$$(1.2) \qquad Y = AX + \varepsilon,$$

where $X$ is the unknown parameter vector from some Hilbert space $\mathcal{H}_X$, $A: \mathcal{H}_X \to \mathcal{H}_Y$ is a linear operator and $\varepsilon$ is a random Gaussian noise in $\mathcal{H}_Y$ with the known correlation structure given by the covariance operator $\Sigma$. The goal is to estimate a linear functional $\theta = \theta(X)$ that can be represented in the form $\langle \vartheta, X \rangle$ for some known element $\vartheta \in \mathcal{H}_X$. Such problems are usually considered as more complex than the usual nonparametric regression estimation due to the poor rate of estimation. Moreover, the difficulty that is usually associated with the attained estimation accuracy increases with the degree of illposedness. A naive estimation approach is based on the explicit least-square solution of the problem (1.2):

$$\widetilde{\theta} = \langle \vartheta, (A^*A)^- A^* Y \rangle = \langle A(A^*A)^- \vartheta, Y \rangle = \langle \phi, Y \rangle,$$

where $A^*$ is the conjugate operator to $A$, $C^-$ means a pseudo-inverse of $C$ and $\phi = A(A^*A)^-\vartheta$. However, this approach cannot be efficiently applied if $A$ is a compact operator because the inverse of $A^*A$ does not exist or is an unbounded operator. One can regularize the problem if some additional information about smoothness of the element $X$ is available. This allows to replace $(A^*A)^-$ by its regularization $g_\alpha(A^*A)$ where $g_\alpha$ means some regularized inversion and $\alpha$ is the corresponding parameter. See, for example, Goldenshluger and Pereversev (2003), Goldenshluger (1999) and Goldenshluger and Pereversev (2000) for typical examples. The quality of estimation heavily depends on the choice of the regularization parameter $\alpha$ and its choice is a challenging problem. Usually one fixes a finite ordered set of values $\alpha_1 < \alpha_2 < \cdots < \alpha_K$ and considers the corresponding estimates

$$\widetilde{\theta}_k = \langle \phi_k, Y \rangle, \qquad \phi_k = Ag_{\alpha_k}(A^*A)\vartheta.$$

Now the original problem can be reformulated as follows: given a set of estimates $\widetilde{\theta}_k$ for known vectors $\phi_k$, build an estimate $\widehat{\theta}$ of the functional $\theta$ that performs nearly as good as the best in this family. We present one particular example for the considered set-up borrowed from Golubev (2004). More examples include a positron emission tomography problem, Cavalier (2001), functional data analysis, Cai and Hall (2006), among many others.

Our analysis focuses on demonstrating the oracle efficiency of the constructed adaptive procedure rather than on establishing the optimal rate of convergence on functional classes. The mentioned efficiency of any adaptive (data-driven) method can be measured by the ratio of the risk of the proposed method to the "oracle" risk which corresponds to the optimal choice of the regularization parameter for the model at hand. One message of this note is that this statistical part of the linear inverse problem is actually not harder than in the classical nonparametric inference. Moreover, in the inverse problem set-up it is typically easier to do a statistical adaption because the likelihood profile is not so flat as in the classical nonparametric regression. In some examples presented in our simulation study in Section 4, the risk of the adaptive procedure is even smaller than the oracle risk.

1.3. *Example for a sequence space model.* We consider the statistical problem with observations $y_1, \ldots, y_M$ following the "sequence space" equation

(1.3) $$y_i = \mu_i + \sigma_i \varepsilon_i, \qquad i = 1, \ldots, M,$$

where $\varepsilon_i$ are independent standard normal and the standard deviations $\sigma_i$ are known while the mean values $\mu_i$ are unknown. The variances $\sigma_i^2$ are usually constant for the regression set-up or grow with $i$ for ill-posed inverse problems.

One particular problem in this set-up can be to estimate the sum

$$\theta = \mu_1 + \cdots + \mu_M,$$

where $M$ can be equal to infinity assuming that the sum of the $\mu_i$'s is absolutely convergent. The "naive" estimate $\widetilde{\theta} = \sum_{i=1}^{M} y_i$, even for a finite $M$, has a very large variance $\sum_{i=1}^{M} \sigma_i^2$ and hence, can be highly inefficient. The smoothing idea leads to the set of the spectral cut-off estimates

$$\widetilde{\theta}_k = \langle \phi_k, Y \rangle = y_1 + \cdots + y_{m_k},$$

where $\phi_k = (1, \ldots, 1, 0, \ldots, 0)$ is the vector with the first $m_k$ entries equal to one and the others equal to zero, while $m_k$ is a fixed decreasing sequence of finite indices $M \geq m_1 > m_2 > \cdots > m_K \geq 1$.

One can easily compute for $k = 1, \ldots, K$ and $l < k$

$$\theta_k \stackrel{\text{def}}{=} \mathbf{E}\widetilde{\theta}_k = \mu_1 + \cdots + \mu_{m_k}, \qquad \xi_k \stackrel{\text{def}}{=} \widetilde{\theta}_k - \theta_k = \sigma_1 \varepsilon_1 + \cdots + \sigma_{m_k} \varepsilon_{m_k}$$

and $v_k \stackrel{\text{def}}{=} \mathrm{Var}(\widetilde{\theta}_k) = \sigma_1^2 + \cdots + \sigma_{m_k}^2$. The major difficulty in applying the smoothing approach is the proper choice of the parameter $k$ or, equivalently the cutting point $m_k$. Small values of $k$ lead to a huge variance $v_k$ of the estimate $\widetilde{\theta}_k$ while large $k$-values can result in a big bias $b_k = \theta - \theta_k = \sum_{i=m_k+1}^{M} \mu_i$. The "oracle" choice balances the approximation and stochastic errors. However, this ideal choice assumes that the bias (the approximation error) is known. The problem we consider in this paper is to develop an adaptive (data-driven) choice that mimics the "oracle" and achieves the best possible performance among the set of estimates $\widetilde{\theta}_k$.

**2. Description of the method.** This section presents the considered adaptive estimation procedure. We first describe some simple properties of the estimates $\widetilde{\theta}_k$ that will be used in the construction. Then we present the adaptive estimation method.

The definition $\widetilde{\theta}_k = \theta_k + \xi_k$ for any $k \leq K$ yields $\mathbf{E}\widetilde{\theta}_k = \theta_k$ and $\mathrm{Var}\,\widetilde{\theta}_k = \mathbf{E}\xi_k^2 = v_k$. Moreover, $\xi_k$ is a zero-mean Gaussian random variable and for any $r > 0$ and any $\lambda < 1$

(2.1) $$\mathbf{E}|v_k^{-1}(\widetilde{\theta}_k - \theta_k)^2|^r = \mathfrak{c}_r,$$

(2.2) $$\mathbf{E}\exp\{\lambda v_k^{-1}(\widetilde{\theta}_k - \theta_k)^2/2\} = (1-\lambda)^{-1/2},$$

where $\mathfrak{c}_r = \mathbf{E}|\xi|^{2r}$ and $\xi$ is standard normal. Due to this result, $\widetilde{\theta}_k$ is a reasonable estimate of $\theta$ if the bias $\theta_k - \theta$ is sufficiently small relative to the standard deviation $v_k^{1/2}$. In particular, in the "no bias" situation $\theta_k = \theta$ the estimate $\widetilde{\theta}_k$ leads to the accuracy of order $v_k^{1/2}$ and one can build confidence intervals for the parameter $\theta_k$ in the form

(2.3) $$\mathcal{E}_k(\mathfrak{z}) = \{u : v_k^{-1}(\widetilde{\theta}_k - u)^2/2 \leq \mathfrak{z}\}.$$

If $\mathfrak{z}$ is sufficiently large, then the result (2.2) ensures that $\mathcal{E}_k(\mathfrak{z})$ contains $\theta_k$ with a high probability.

2.1. *Adaptive choice of an estimate out of a given family.* Our starting point is the given family of estimates $\widetilde{\theta}_k$ for $k = 1, \ldots, K$ ordered by their variability so that the variance $v_k$ of $\widetilde{\theta}_k$ decreases with $k$. We aim to select a data-driven index $\widehat{k}$ or equivalently the estimate $\widehat{\theta} = \widetilde{\theta}_{\widehat{k}}$, which minimizes the corresponding estimation risk.

For a given sequence of estimates $\widetilde{\theta}_k = \theta_k + \xi_k$ consider the sequence of nested hypothesis $H_k: \theta_1 = \cdots = \theta_k = \theta$. The procedure is sequential: we start with $k = 2$ and at every step $k$ the hypothesis $H_k$ is tested against $H_1, \ldots, H_{k-1}$. If $H_k$ is not rejected then we continue with the next larger $k$. The final estimate corresponds to the latest accepted hypothesis. For testing $H_k$ against $H_l$ with $l < k$, we check that the new estimate $\widetilde{\theta}_k$ belongs to the confidence intervals built on the base of $\widetilde{\theta}_l$. More precisely, we apply the test statistics:

$$T_{lk} = v_l^{-1}(\widetilde{\theta}_l - \widetilde{\theta}_k)^2/2, \qquad l < k,$$

where $v_l$ is the variance of $\widetilde{\theta}_l$. Big values of $T_{lk}$ indicate a significant difference between the estimates $\widetilde{\theta}_l$ and $\widetilde{\theta}_k$. Due to the definition (2.3), the event $A_{lk} = \{T_{lk} \leq \mathfrak{z}_l\}$ means that $\widetilde{\theta}_k$ belongs to the confidence set $\mathcal{E}_l(\mathfrak{z}_l)$ based on $\widetilde{\theta}_l$. The estimate $\widetilde{\theta}_k$ (or the hypothesis $H_k$) is accepted if $H_{k-1}$ was accepted and $T_{lk} \leq \mathfrak{z}_l$ for all $l < k$, that is, the new estimate $\widetilde{\theta}_k$ belongs to the intersection of all the confidence intervals $\mathcal{E}_l(\mathfrak{z}_l)$ built on the previous steps of the procedure. The formal definition is given by

$$\widehat{k} = \max\{k \leq K : T_{lm} \leq \mathfrak{z}_l, \ \forall l < m \leq k\}.$$

Here the "critical values" $\mathfrak{z}_1, \ldots, \mathfrak{z}_{K-1}$ are the parameters of the procedure. Their choice is discussed in Section 2.2.

The selected random index $\widehat{k}$ means the largest accepted $k$. The corresponding adaptive estimate $\widehat{\theta}$ is $\widetilde{\theta}_{\widehat{k}} : \widehat{\theta} = \widetilde{\theta}_{\widehat{k}}$. We also define the adaptive estimate $\widehat{\theta}_k$ as the latest accepted after the first $k$ steps:

$$\widehat{\theta}_k = \widetilde{\theta}_{\min\{\widehat{k}, k\}}.$$

The described procedure involves $K - 1$ parameters and their automatic choice is ultimately required for practical applications of the method. Our next step is the method for an automatic selection of the critical values $\mathfrak{z}_k$.

2.2. *Choice of the critical values $\mathfrak{z}_k$ using a "propagation condition."* The way of selecting the critical values $\mathfrak{z}_1, \ldots, \mathfrak{z}_{K-1}$ is similar to the standard approach of hypothesis testing theory: to provide the prescribed performance of the procedure under the simplest (null) hypothesis. In the considered set-up, the null means $\theta_1 = \theta_2 = \cdots = \theta_K = \theta$. We will show below in Theorem 3.3 that the particular value of $\theta$ is unimportant and it suffices to only consider $\theta = 0$. In what follows we denote by $\mathbf{P}_0$ the distribution of the data in this situation and $\mathbf{E}_0$ means the corresponding mathematical expectation. By the definition of the procedure, accepting $H_k$ for some $k \leq K$ yields $\widehat{\theta}_k = \widetilde{\theta}_k$ and rejecting of $H_k$ means $\widehat{\theta}_k \neq \widetilde{\theta}_k$. We refer to the latter as a "false alarm" because the procedure terminates in the situation where it should not. If such false alarms occur too often, it is an indication that the critical values $\mathfrak{z}_k$ are not large enough. The usual $\alpha$-level condition on any testing procedure is that under the null it rejects the null hypothesis with the probability not exceeding $\alpha$. For the considered multiple test procedure this condition reads as $\mathbf{P}_0(\widehat{\theta}_k \neq \widetilde{\theta}_k) \leq \alpha$. We slightly modify this condition to adapt it to the problem of adaptive estimation by selecting a polynomial loss function instead of the indicator of the error decision. Rejecting the null hypothesis happens if $v_l^{-1}(\widetilde{\theta}_l - \widetilde{\theta}_k)^2/2 > \mathfrak{z}_l$, in which can be interpreted that the estimate $\widetilde{\theta}_k$ does not belong to the confidence interval $\mathcal{E}_l(\mathfrak{z}_l)$ built on the base of $\widetilde{\theta}_l$. In the testing problem it only matters how often such false alarms occur. In the considered problem of adaptive estimation we focus on the risk associated with such a false alarm. Therefore, the particular indices $l, k$ and the size of $v_l^{-1}(\widetilde{\theta}_l - \widetilde{\theta}_k)^2$ matter as well. Suppose that some loss power $r > 0$ is fixed. By (2.1)

$$\mathbf{E}_0 |v_l^{-1}(\widetilde{\theta}_l - \theta)^2|^r = \mathfrak{c}_r, \qquad \text{for all } l \leq K,$$

where $\mathfrak{c}_r = E|\xi|^{2r}$ and $\xi$ is standard normal. We require that the parameters $\mathfrak{z}_1, \ldots, \mathfrak{z}_{K-1}$ of the procedure are selected in such a way that

$$(2.4) \qquad \mathbf{E}_0 |v_k^{-1}(\widehat{\theta}_k - \widetilde{\theta}_k)^2|^r \leq \alpha \mathfrak{c}_r, \qquad k = 2, \ldots, K.$$

The meaning of this condition is that at every step $k$ of the procedure, the risk associated with false alarms is at most an $\alpha$-fraction of the best-possible estimation risk. Here $\alpha$ is the preselected constant, which is similar to the confidence level of a testing procedure. This gives us $K - 1$ conditions to fix $K - 1$ parameters. As in the testing problem, we are interested to select the critical values as small as possible under the constraint (2.4). Note that the choice $r$ very close to zero leads back to the indicator loss function $\mathbf{1}(\widehat{\theta}_k \neq \widetilde{\theta}_k)$ and thus, to the usual error of the first kind for the multiple testing procedure.

Our definition still involves two parameters $\alpha$ and $r$. It is important to mention that their choice is subjective and there is no way for an automatic selection in the considered local or pointwise set-up. Moreover, the possibility of tuning such parameters in particular applications is an important advantage of the approach. Our aim is to develop a procedure that combines and balances two important features: stability in the parametric situation and sensitivity under deviations from the parametric null hypothesis. The propagation condition (2.4) is exactly a constraint on the stability in the parametric case, and we aim to optimize the sensitivity of the method under this constraint. A proper choice of the power $r$ for the loss function as well as the "confidence level" $\alpha$ depends on the particular application and on the additional subjective requirements to the procedure. Taking a large $r$ and small $\alpha$ would result in an increase of the critical values and therefore improves the performance of the method in the parametric situation at cost of some loss of sensitivity to deviations from the parametric situation. This behavior is analogous to the hypothesis testing problem where a small $\alpha$ reduces the first-kind error at costs of the test's power. Theorem 3.1 presents some upper bounds for the critical values $\mathfrak{z}_k$ as functions of $\alpha$ and $r$ in the form $a_1 \log(K) + a_2\{\log \alpha^{-1} + r(K-k)\}$ with some coefficients $a_1$ and $a_2$. We see that these bounds linearly depend on $r$ and on $\log \alpha^{-1}$. For our examples, we apply a relatively small value $r = 1/2$. We also apply $\alpha = 1$ although the other values in the range $[0.5, 1]$ lead to very similar results. It is worth mentioning that both the procedure and the theoretical study apply and lead to reasonable results whatever $r$ and $\alpha$ are. This makes a striking difference with many other proposals; see the references in the introduction for selecting the tuning parameter(s). Typically one requires that the critical values (thresholds) $\mathfrak{z}$ are sufficiently large and the theory is only valid under this condition.

2.3. *Sequential choice.* The set of conditions (2.4) does not directly define the critical values $\mathfrak{z}_k$. We present below one sequential method for fixing $\mathfrak{z}_k$ one after another starting from $\mathfrak{z}_1$. The idea is to provide that the relative impact of each $\mathfrak{z}_k$ in the total risk in (2.4) is the same for every $k \leq K-1$. We start with $\mathfrak{z}_1$ and set $\mathfrak{z}_2 = \cdots = \mathfrak{z}_{K-1} = \infty$. This effectively means that every new estimate $\widetilde{\theta}_k$ is only compared with $\widetilde{\theta}_1$. We run the procedure with such critical values. The resulting adaptive estimate after step $k$ is denoted by $\widehat{\theta}_k(\mathfrak{z}_1)$. We select $\mathfrak{z}_1$ as the minimal value providing

$$(2.5) \qquad \mathbf{E}_0 |v_k^{-1}\{\widehat{\theta}_k(\mathfrak{z}_1) - \widetilde{\theta}_k\}^2|^r \leq \frac{1}{K-1}\alpha \mathfrak{c}_r, \qquad k = 2,\ldots,K.$$

Such a value exists because the choice $\mathfrak{z}_1 = \infty$ leads to $\widehat{\theta}_k = \widetilde{\theta}_k$ for all $k$.

Similarly, we specify $\mathfrak{z}_2$ by considering the situation with the previously fixed $\mathfrak{z}_1$, some finite $\mathfrak{z}_2$ and all the remaining critical values equal to infinity, and so on. For the formal definition, suppose that $\mathfrak{z}_1,\ldots,\mathfrak{z}_{m-1}$ have been already fixed for some $m > 1$ and define for any $\mathfrak{z}_m$ the adaptive estimates $\widehat{\theta}_k(\mathfrak{z}_1,\ldots,\mathfrak{z}_m)$ for $k > m$, which come out of the procedure with the critical values $(\mathfrak{z}_1,\ldots,\mathfrak{z}_m,\infty,\ldots,\infty)$. We select $\mathfrak{z}_m$ as the minimal value providing

$$(2.6) \qquad \mathbf{E}_0 |v_k^{-1}\{\widehat{\theta}_k(\mathfrak{z}_1,\ldots,\mathfrak{z}_m) - \widetilde{\theta}_k\}^2|^r \leq \frac{m}{K-1}\alpha \mathfrak{c}_r, \qquad k = m+1,\ldots,K.$$

Such a value exists because the choice $\mathfrak{z}_m = \infty$ leads to $\widehat{\theta}_k(\mathfrak{z}_1, \ldots, \mathfrak{z}_m) = \widehat{\theta}_k(\mathfrak{z}_1, \ldots, \mathfrak{z}_{m-1})$ and even a stronger condition has been already checked at the previous step.

The condition (2.5) describes the impact of the first critical value in the risk (2.4) while (2.6) describes the accumulated impact of the first $m$ critical values. The factor $m/(K-1)$ in the right-hand side of (2.6) is chosen to ensure that every critical value $\mathfrak{z}_k$ has the same impact.

Our construction guarantees that the selected sequence $\mathfrak{z}_k$ is minimal under the set of conditions (2.6) in the sense that one cannot select another sequence $\mathfrak{z}'_k < \mathfrak{z}_k$ for all $k$ such that (2.6) is still fulfilled. Indeed, let $\{\mathfrak{z}'_k\}$ be another sequence that ensures (2.6) and let $m$ be the first index for which $\mathfrak{z}'_m < \mathfrak{z}_m$. Then the condition

$$\mathbf{E}_0 |v_k^{-1}\{\widehat{\theta}_k(\mathfrak{z}'_1, \ldots, \mathfrak{z}'_{m-1}, \mathfrak{z}_m) - \widetilde{\theta}_k\}^2|^r \leq \frac{m}{K-1} \alpha \mathfrak{c}_r, \qquad k > m,$$

on $\mathfrak{z}_m$ is even stronger than (2.6) and one cannot select $\mathfrak{z}'_m < \mathfrak{z}_m$ to ensure it. This contradiction shows minimality of the sequence $\mathfrak{z}_k$.

An explicit form for the critical values $\mathfrak{z}_k$ is not available but they can be easily computed using the Monte Carlo simulations from the null hypothesis; see Section 4 for details.

**3. Theoretical study.** This section presents some properties of the adaptive estimate $\widehat{\theta}$ of the target value $\theta$. We suppose that the parameters $\mathfrak{z}_k$ of the procedure are selected in such a way that the condition (2.4) is fulfilled. The main result is the "oracle" property of the adaptive estimate $\widehat{\theta}$, which claims that the risk of adaptive estimation is up to some multiplier as good as the risk of the ideal (oracle) estimate. This multiplier is directly related to the applied critical values $\mathfrak{z}_k$ and in typical situations it is at most logarithmic in the sample size. In the proof we distinguish between three cases: parametric, local parametric and nonparametric. The parametric case means that $\theta_k \stackrel{\text{def}}{=} \mathbf{E}\widetilde{\theta}_k \equiv \theta$ for all $k \leq K$. This case can be easily reduced to the null hypothesis $\theta_1 = \cdots = \theta_K = 0$ and the oracle property of the adaptive estimate $\widehat{\theta}$ is ensured by the construction, more precisely, by the propagation condition (2.4). The local parametric case means that for some $k < K$ holds $\theta_1 = \cdots = \theta_k = \theta$. In this case, the construction ensures the oracle property for the adaptive estimate $\widehat{\theta}_k$ obtained after the first $k$ steps of the procedure. Then we show that a similar oracle property of the estimate $\widehat{\theta}_k$ can be obtained in the nonparametric situation under the so-called "small modeling bias" condition. This condition is used to give a formal definition of the oracle choice. The final oracle result for the adaptive estimate $\widehat{\theta}$ is obtained by combining the previously established "propagation" result under the small modeling bias condition with the "stability" property, which is ensured by the adaptive procedure itself.

3.1. *Bounds for the critical values.* This section presents some upper and lower bounds for the critical values $\mathfrak{z}_k$. The results are established under the following condition on the variances $v_k$.

(*MD*) for some constants $\mathfrak{u}_0, \mathfrak{u}$ with $1 < \mathfrak{u}_0 \leq \mathfrak{u}$, the variances $v_k$ satisfy

$$v_{k-1} \leq \mathfrak{u} v_k, \qquad \mathfrak{u}_0 v_k \leq v_{k-1}, \qquad 2 \leq k \leq K.$$

We also denote for $l < k \leq K$

$$v_{l,k} \stackrel{\text{def}}{=} \operatorname{Var}(\widetilde{\theta}_l - \widetilde{\theta}_k).$$

Our first result presents some upper bound for the parameters $\mathfrak{z}_k$ under condition (*MD*). The proof is given in the Appendix.

THEOREM 3.1. *Assume* (*MD*). *Let $\gamma$ be such that for all $l < k \leq K$*

$$v_{l,k}/v_l \leq \gamma. \tag{3.1}$$

*Then there is a constant $C_1$ depending on $r$, $\mathfrak{u}_0$, and $\mathfrak{u}$ only such that the choice*

$$\mathfrak{z}_k = \gamma \{\log \alpha^{-1} + r \log(v_k/v_K)\} + C_1 \log K$$

*ensures (2.4) for all $k \leq K$. Particularly, $\mathbf{E}_0 |v_K^{-1}(\widetilde{\theta}_K - \widehat{\theta})^2|^r \leq \alpha \mathfrak{c}_r$.*

REMARK 3.1. The result of Theorem 3.1 presents some upper bounds for the critical values. These upper bounds will be used for our theoretical study; however, they do not appear in the proposed adaptive procedure. An interesting observation is that these upper bounds linearly decrease with $k$. Indeed, by condition (*MD*) $\log(v_k/v_K) \leq (K-k) \log \mathfrak{u}$ and $\log(v_k/v_K) \geq (K-k) \log \mathfrak{u}_0$. The reason for a decrease of $\mathfrak{z}_k$ with $k$ can be explained as follows. Under the null hypothesis the procedure should not terminate at intermediate steps and the oracle estimate is $\widetilde{\theta}_K$. An early stop ("false alarm") $\widehat{k} = k$ for $k < K$ results in selecting the estimate $\widetilde{\theta}_k$, which has much larger variability than $\widetilde{\theta}_K$. The smaller $k$ is, the larger is the associated loss in the estimation quality. Therefore, the test at the early stage of the procedure should be rather conservative while a "false alarm" at the final steps of the procedure is not so critical, and we are more interested to improve sensitivity by applying nonconservative critical values.

Our next result shows that the linear growth of the critical values $\mathfrak{z}_k$ with $K - k$ is not only sufficient but also necessary for providing (2.4). To highlight the contribution of every particular value $\mathfrak{z}_k$, we consider the situation when all the previous parameters are equal to infinity: $\mathfrak{z}_1 = \cdots = \mathfrak{z}_{k-1}$. This effectively means that the procedure cannot terminate at the first $k-1$ steps due to a possibly wrong choice of the corresponding critical values.

THEOREM 3.2. *Assume* (MD). *Suppose that for a fixed $k < K$, it holds $\mathfrak{z}_1 = \cdots = \mathfrak{z}_{k-1} = \infty$. Then the condition (2.4) implies that*

$$\mathfrak{z}_k \geq \frac{v_{k,k+1}}{v_k}\{r\log(v_{k,K}/v_K) + \log\alpha^{-1} - C_2\log(K)\}$$

*for some positive constant $C_2$ depending on $r, \mathfrak{u}, \mathfrak{u}_0$ only.*

The proof is again moved to the Appendix.

REMARK 3.2. Our main oracle result particularly shows that the leading term in the risk linearly depends on the value $\mathfrak{z}_{k^*}$ where $k^*$ is the optimally selected index. Therefore, obtaining a sharp oracle result would require bringing together the upper and lower bounds for the critical values $\mathfrak{z}_k$. In our results these two bounds differ by the factor $\gamma v_k/v_{k,k+1}$ with $\gamma$ from (3.1). The value $\gamma$ is usually close to one because the estimates $\widetilde{\theta}_k$ are positively correlated with each other in the most of cases. However, the value $v_{k,k+1}/v_k$ can be small by the same reason. So, obtaining a sharp oracle result would require some modification of the presented procedure; cf. Lepski and Spokoiny (1997). The further discussion of this issue lies beyond the scope of this paper.

3.2. *Behavior in the local parametric situation.* The parametric situation can be understood as the case when $\theta_1 = \theta_2 = \cdots = \theta_K$. In this case the estimate $\widetilde{\theta}_K$ is unbiased and has the smallest variance and hence, the smallest risk described by the formula $\mathbf{E}|v_K^{-1}(\widetilde{\theta}_K - \theta)^2|^r = \mathfrak{c}_r$. A natural requirement to any adaptive procedure is to provide a similar accuracy of the adaptive estimate under the parametric hypothesis. Similarly, the local parametric situation corresponds to the case when $\theta_1 = \cdots = \theta_k = \theta$ for some $k \leq K$. In this case it is natural to require that the adaptive estimate $\widehat{\theta}_k$ after $k$ steps is close to its nonadaptive counterpart $\widetilde{\theta}_k$. This property is actually provided by the construction of the critical values.

THEOREM 3.3. *Let $\theta_1 = \theta_2 = \cdots = \theta_K = \theta$. Then it holds*

$$\mathbf{E}|v_K^{-1}(\widehat{\theta} - \widetilde{\theta}_K)^2|^r \leq \alpha\mathfrak{c}_r.$$

*Moreover, if $\theta_1 = \theta_2 = \cdots = \theta_k = \theta$ for some $k \leq K$, then*

$$\mathbf{E}|v_k^{-1}(\widehat{\theta}_k - \widetilde{\theta}_k)^2|^r \leq \alpha\mathfrak{c}_r.$$

PROOF. Only the differences $\widetilde{\theta}_l - \widetilde{\theta}_k$ appear in the definition of the test statistics $T_{lk}$. In view of the decomposition $\widetilde{\theta}_k = \theta + \xi_k$, the value $\theta$ cancels there. Similarly, the adaptive estimate $\widehat{\theta}_k$ coincides with one of $\widetilde{\theta}_1, \ldots, \widetilde{\theta}_k$ and the value $\theta$ cancels in the difference $\widehat{\theta}_k - \widetilde{\theta}_k$ as well. Hence, we can assume $\theta = 0$ and $\widetilde{\theta}_k = \xi_k$. Then the results follow from the constraints (2.4) on the critical values $\mathfrak{z}_k$. □

3.3. *"Small modeling bias" condition and "propagation" property.* Theorem 3.3 describes the performance of the estimate $\widetilde{\theta}_k$ under the parametric or local parametric assumption. Now we aim to extend this result to the general nonparametric situation when the identities $\theta_1 = \theta_2 = \cdots = \theta_k = \theta$ are only approximately fulfilled and the deviation from the null hypothesis $H_k$ is not significant.

As mentioned in Section 2.2, the choice of critical values $\mathfrak{z}_k$ is determined by the joint distribution of the test statistics $T_{lk} = v_l^{-1}(\widetilde{\theta}_l - \widetilde{\theta}_k)^2$ under the measure $\mathbf{P}_0$ corresponding to the parametric hypothesis $\theta_1 = \theta_2 = \cdots = \theta_K = 0$. An extension of this result to the nonparametric situation leads to considering the similar distribution in the general case. Let $\mathbf{P}_k$ mean the joint distribution of $\widetilde{\theta}(k) = (\widetilde{\theta}_1, \ldots, \widetilde{\theta}_k)^\top$ for $k \geq 1$. By the model assumption, this is a Gaussian vector with $\mathbf{E}\widetilde{\theta}(k) = \theta(k) = (\theta_1, \ldots, \theta_k)^\top$. Let also $B_k$ be the covariance matrix of the vector $\widetilde{\theta}(k)$. Then $\mathbf{P}_k$ is the normal distribution with the mean $\theta(k)$ and the covariance matrix $B_k$, $\mathbf{P}_k = \mathcal{N}(\theta(k), B_k)$. Similarly, $\mathbf{P}_{\theta,k}$ denotes the distribution of $\widetilde{\theta}(k)$ under the local parametric situation $\theta_1 = \cdots = \theta_k = \theta$, that is, $\mathbf{P}_{\theta,k} = \mathcal{N}(\theta_0(k), B_k)$, where $\theta_0(k) = (\theta, \ldots, \theta)^\top$. Let $b(k) = (b_1, \ldots, b_k)^\top$ with $b_k = \theta_k - \theta$.

LEMMA 3.1. *For $k \geq 1$, define*
$$\Delta_k \stackrel{\text{def}}{=} b^\top(k) B_k^{-1} b(k).$$
*Then the Kullback–Leibler divergence $\mathcal{K}(\mathbf{P}_k, \mathbf{P}_{\theta,k})$ fulfills*
$$\mathcal{K}(\mathbf{P}_k, \mathbf{P}_{\theta,k}) \stackrel{\text{def}}{=} \mathbf{E}_k \log\left(\frac{d\mathbf{P}_k}{d\mathbf{P}_{\theta,k}}\right) = \Delta_k/2$$
*and the values $\Delta_k$ grow with $k$. It also holds for any $s > 1$*
$$\frac{1}{s} \log \mathbf{E}_{\theta,k}\left(\frac{d\mathbf{P}_k}{d\mathbf{P}_{\theta,k}}\right)^s = \frac{\Delta_k(s-1)}{2}.$$
*Moreover, if $\zeta$ is measurable function of $\widetilde{\theta}_1, \ldots, \widetilde{\theta}_k$, then with $s' = s/(s-1)$*
$$\mathbf{E}\zeta \leq (\mathbf{E}_{\theta,k}\zeta^{s'})^{1/s'} \exp\{\Delta_k(s-1)/2\}.$$
*In particular, for $s = 2$ it holds $\mathbf{E}\zeta \leq (\mathrm{e}^{\Delta_k} \mathbf{E}_{\theta,k}\zeta^2)^{1/2}$.*

PROOF. Define $Z_k = d\mathbf{P}_k/d\mathbf{P}_{\theta,k}$. Then
$$\log Z_k = b^\top(k) B_k^{-1/2} \xi_k + b^\top(k) B_k^{-1} b(k)/2$$
with $\xi_k \sim \mathcal{N}(0,1)$ and hence $\mathbf{E}_k \log(Z_k) = \Delta_k/2$. Therefore, $\Delta_k$ is twice the Kullback–Leibler divergence between two measures $\mathbf{P}_k$ and $\mathbf{P}_{\theta,k}$ obtained by projecting the measures $\mathbf{P}$ and $\mathbf{P}_\theta$ on the $\sigma$-field generated by $\widetilde{\theta}_1, \ldots, \widetilde{\theta}_k$ and growing with $k$. This immediately implies that $\Delta_k$ monotonously increase with $k$, that is, $\Delta_k \leq \Delta_{k'}$ for $k < k'$. Similarly,

$$\begin{aligned}\mathbf{E}_{\theta,k} Z_k^s &= \mathbf{E}_{\theta,k} \exp\{sb^\top(k) B_k^{-1/2} \xi_k - b^\top(k) B_k^{-1} b(k) s/2\} \\ &= \exp\{b^\top(k) B_k^{-1} b(k) (s^2 - s)/2\}.\end{aligned}$$

Next, let $\zeta$ be a measurable function of the vector $\widetilde{\theta}(k)$. It holds $\mathbf{E}\zeta = \mathbf{E}_{\theta,k} \zeta Z_k$. By the Hölder inequality

$$\mathbf{E}_{\theta,k} \zeta Z_k \leq (\mathbf{E}_{\theta,k} \zeta^{s'})^{1/s'} (\mathbf{E}_{\theta,k} Z_k^s)^{1/s}$$

and the assertion follows. □

Due to Lemma 3.1, the value $\Delta_k$ can be used to measure the distance between the two models: one corresponds to the local parametric situation with $\theta_1 = \theta_2 = \cdots = \theta_k = \theta$ and the other describes the distribution of the same vector $\widetilde{\theta}(k)$ in the general nonparametric situation. We call this value $\Delta_k$ *the modeling bias* because it describes how much we have to pay in the risk for using the "wrong" parametric model in place of the underlying nonparametric one. The "small modeling bias" (SMB) condition simply means that the value $\Delta_k$ does not exceed some sufficiently small value $\Delta$.

The result of Lemma 3.1 implies that the bound for the risk of estimation $\mathbf{E}_0\{v_k^{-1}(\widetilde{\theta}_k - \theta)^2\}^r$ under the parametric hypothesis translates under the SMB condition $\Delta_k \leq \Delta$ into the bound for the risk $\mathbf{E}\{v_k^{-1}(\widetilde{\theta}_k - \theta)^2\}^{r/s'}$. Similarly one can bound $\mathbf{E}\{v_k^{-1}(\widehat{\theta}_k - \widetilde{\theta}_k)^2\}^{r/s'}$.

In what follows we apply the result of Lemma 3.1 with $s = s' = 2$, which nicely simplifies the notation. Note, however, that any $s > 1$ can be used. For instance, taking a large $s$ leads to the value of $s'$ close to one.

THEOREM 3.4. *For any $r > 0$, it holds for every $k \leq K$*

$$\mathbf{E}\{v_k^{-1}(\widetilde{\theta}_k - \theta)^2\}^{r/2} \leq \sqrt{e^{\Delta_k} \mathfrak{c}_r},$$

$$\mathbf{E}\{v_k^{-1}(\widetilde{\theta}_k - \widehat{\theta}_k)^2\}^{r/2} \leq \sqrt{e^{\Delta_k} \alpha \mathfrak{c}_r}.$$

The bound follows directly from Lemma 3.1 and Theorem 3.3.

We call this result the "propagation" property because it ensures that with a high probability the procedure does not terminate yielding $\widehat{\theta}_k = \widetilde{\theta}_k$ as long as the SMB condition $\Delta_k \leq \Delta$ is fulfilled. Note that a similar property has been proved for the original procedure in Lepski (1990) and Lepski (1991, 1992), however, under the additional condition that the critical values $\mathfrak{z}_k$ are sufficiently large. We instead use the propagation condition (2.4) and the SMB condition.

3.4. *"Stability after propagation" and oracle results.* Due to the "propagation" result of Theorem 3.4, the procedure performs well as long as the SMB condition is fulfilled, which means that the value $\Delta_k$ remains bounded by some (small) constant. We formalize this condition in the form $\Delta_k \leq \Delta$. Here $\Delta$ is an arbitrary number that will determine the oracle choice. We will show in Section 3.6 that in typical situations this value $\Delta$ is similar to the ratio of the squared bias to the variance of $\widetilde{\theta}_k$. Note however, that the value $\Delta$ only appears in our theoretical study; it does not affect the procedure. The results apply whatever $\Delta > 0$.

To establish the accuracy result for the final estimate $\widehat{\theta}$, we have to check that the adaptive estimate $\widehat{\theta}_k$ does not vary much at the steps at which the modeling bias $\Delta_k$ becomes large.

THEOREM 3.5 (Stability). *It holds for every $k < K$*

$$(3.2) \qquad v_k^{-1}(\widetilde{\theta}_k - \widehat{\theta})^2 \mathbf{1}(\widehat{k} > k) \leq 2\mathfrak{z}_k.$$

PROOF. The result follows by the definition of $\widehat{\theta} = \widetilde{\theta}_{\widehat{k}}$ and $\widehat{\theta}_k = \widetilde{\theta}_{\min\{\widehat{k},k\}}$ because $\widehat{k}$ is accepted and $\min\{\widehat{k}, k\} \leq \widehat{k}$. □

Combination of the "propagation" and "stability" statements implies the main result concerning the properties of the adaptive estimate $\widehat{\theta}$. In the formulation of this and the further results we assume some constant $\Delta > 0$ to be fixed. We also assume that our set-up is reasonable in the sense that for the very first model the SMB condition $\Delta_1 \leq \Delta$, or equivalently, $b_1^2 \leq \Delta v_1$, is fulfilled. This enables us to correctly define the ideal index $k^*$.

THEOREM 3.6. *Let $k^*$ be the maximal value $k$ such that $\Delta_k \leq \Delta$. Then*

$$(3.3) \qquad \mathbf{E}|v_{k^*}^{-1}(\widetilde{\theta}_{k^*} - \widehat{\theta})^2|^{r/2} \leq \sqrt{\alpha \mathfrak{c}_r e^\Delta} + (2\mathfrak{z}_{k^*})^{r/2}.$$

PROOF. The events $\mathbf{1}(\widehat{k} > k^*)$ and $\mathbf{1}(\widehat{k} \leq k^*)$ do not overlap and $\widehat{\theta} = \widehat{\theta}_{k^*}$ for $\widehat{k} \leq k^*$. This yields the representation
$$\mathbf{E}|v_{k^*}^{-1}(\widetilde{\theta}_{k^*} - \widehat{\theta})^2|^{r/2} = \mathbf{E}|v_{k^*}^{-1}(\widetilde{\theta}_{k^*} - \widehat{\theta})^2|^{r/2}\mathbf{1}(\widehat{k} > k^*) + \mathbf{E}|v_{k^*}^{-1}(\widetilde{\theta}_{k^*} - \widehat{\theta}_{k^*})^2|^{r/2}.$$
Now the result follows from Theorems 3.4 and 3.5. □

3.5. *Discussion.* Here we discuss some issues related to the stated oracle result.

*"Oracle" quality.* Theorem 3.4 ensures that the estimation loss $v_k^{-1}(\widetilde{\theta}_k - \theta)^2$ is bounded with a high probability if the modeling bias $\Delta_k$ is not too big. The oracle choice $k^*$ is the largest one for which the SMB condition $\Delta_k \leq \Delta$ holds leading to the accuracy $|\widetilde{\theta}_{k^*} - \theta|$ of order $v_{k^*}^{1/2}$. We aim to build an adaptive estimate that delivers the same quality as the oracle one. Theorem 3.6 claims that the difference $\widehat{\theta} - \widetilde{\theta}_{k^*}$ between the adaptive estimate $\widehat{\theta}$ and oracle is indeed of order $v_{k^*}^{1/2}$ up to the factor $\sqrt{2\mathfrak{z}_{k^*}}$.

*The "true" value $\theta$.* The "true" value $\theta$ is not explicitly shown in the oracle inequality (3.3). It only enters in the definition of the modeling bias $\Delta_k$ and thus, in the SMB condition $\Delta_k \leq \Delta$ and in the definition of the oracle choice $k^*$. The oracle bound just compares the optimal choice of the index $k^*$ for the given nonparametric model (1.1) with the adaptive index $\widehat{k}$. In fact, the model (1.1) does not require a "true" value $\theta$ to be defined and the oracle result can be formally applied for any $\theta$. However, in our two basic examples of nonparametric regression and linear function estimation such values are defined in a natural way. The quality of estimation of this value $\theta$ can be easily derived from the oracle bound (3.3). We present the corresponding result about the risk of the adaptive estimate $\widehat{\theta}$ for the special case with $r = 1$. The other values of $r$ can be considered as well, one only has to update the constants depending on $r$. We also assume that $\alpha \leq 1$.

COROLLARY 3.7. *Let $k^*$ be the largest $k$ with $\Delta_k \leq \Delta$. Then*
$$v_{k^*}^{-1/2} \mathbf{E}|\widehat{\theta} - \theta| \leq 2\sqrt{e^\Delta} + \sqrt{2\mathfrak{z}_{k^*}}.$$

PROOF. Just observe that
$$|\widehat{\theta} - \theta| \leq |\widetilde{\theta}_{k^*} - \theta| + |\widetilde{\theta}_{k^*} - \widehat{\theta}|$$
and the result follows from Theorem 3.6 in view of $\mathfrak{c}_1 = 1$. $\square$

*Leading term in the risk.* The risk bound in the presented oracle result consists of two terms. The first one $\sqrt{\alpha \mathfrak{c}_r e^\Delta}$ is just a constant. Moreover, by choosing a small $\alpha$, one can make this term arbitrary small. The other term $(2\mathfrak{z}_{k^*})^{r/2}$ is by the bound of Theorem 3.1 of order $\log K$ and thus, under the assumption $(MD)$, it is logarithmic in the sample size. This implies that asymptotically, as the sample size increases, the leading term in the risk bound is exactly the value $(2\mathfrak{z}_{k^*})^{r/2}$. This particularly explains why the choice of a possibly small critical values is an important issue.

*Payment for adaptation.* Recall that in the parametric situation, the risk $\mathbf{E}|v_{k^*}^{-1}(\widetilde{\theta}_{k^*} - \theta)^2|$ of $\widetilde{\theta}_{k^*}$ is bounded by $\mathfrak{c}_1 = 1$; see (2.2). In the nonparametric situation, the result is only slightly worse. The risk bound includes $\sqrt{2\mathfrak{z}_{k^*}}$, which can be logarithmic in the sample size. In addition, it bounds the absolute loss $|\widehat{\theta} - \theta|$ instead of squared loss. Finally, there is an additional factor $\sqrt{e^\Delta}$, which accounts for the use of a wrong parametric model instead of the real one.

3.6. *SMB condition versus "bias-variance trade-off."* The standard approach for selecting the optimal index $k$ is based on balancing an upper bound $\overline{b}_k$ for the bias $b_k = \theta_k - \theta$ and the standard deviation $v_k^{1/2}$ of the estimate $\widetilde{\theta}_k$, see for example Lepski, Mammen and Spokoiny (1997) or Goldenshluger (1998) for a related discussion. This section shows that under some additional technical assumptions this approach is nearly equivalent to the SMB condition advocated in this paper.

In addition to $(MD)$ we suppose the following properties of the covariance matrices $B_k = \text{Cov}(\widetilde{\theta}(k))$. Let $B_{k,\text{diag}}$ be the diagonal matrix with the same diagonal entries $v_k$ as for $B_k$. Define also $D_k = B_k^{1/2}$ and $D_{k,\text{diag}} = B_{k,\text{diag}}^{1/2}$. The required conditions reads as follows:

$(Dk)$ It holds for some constant $\mathfrak{s}$ and all $k \leq K$

$$D_k^{-1} \preceq \mathfrak{s} D_{k,\text{diag}}^{-1}.$$

Here the notation $A \preceq B$ for two symmetric matrices $A, B$ means that $|Av| \leq |Bv|$ for any vector $v$. If $B$ is invertible, this is equivalent to saying that the maximal eigenvalue of the matrix $B^{-1}A^2B^{-1}$ is bounded by $\mathfrak{s}^2$.

Condition $(Dk)$ allows to rewrite the SMB condition $|D_k^{-1}b(k)|^2 \leq \Delta$ in the following form:

$$(3.4) \qquad |D_{k,\text{diag}}^{-1}b(k)|^2 \equiv b_1^2/v_1 + \cdots + b_k^2/v_k \leq \Delta/\mathfrak{s}^2.$$

Let $\overline{b}_k$ be a monotonously increasing upper bound for $|b_k|$: $\overline{b}_k = \max_{l \leq k} |b_l|$. For the considered problem of pointwise estimation, the bias-variance trade-off is usually written in the form

$$(3.5) \qquad \overline{b}_{k^*} \leq C_b v_{k^*}^{1/2}$$

for some fixed constant $C_b$; see Lepski, Mammen and Spokoiny (1997). The next result shows that this relation implies the SMB condition (3.4).

THEOREM 3.8. *Suppose* $(MD)$ *and* $(Dk)$. *Then for the index $k^*$ defined by the balance relation (3.5), the SMB condition $\Delta_{k^*} \leq \Delta$ is also fulfilled with* $\Delta = \mathfrak{s}^2 C_{\mathfrak{u}_0} C_b^2$.

PROOF. Let $k$ be such that $\overline{b}_k \leq C_b v_k^{-1/2}$. Then

$$b_1^2/v_1 + \cdots + b_k^2/v_k \leq \overline{b}_k^2(v_1^{-1} + \cdots + v_k^{-1}) \leq C_{\mathfrak{u}_0} \overline{b}_k^2 v_k^{-1},$$

$C_{\mathfrak{u}_0} = (1 - \mathfrak{u}_0^{-1})^{-1}$. Now condition $(Dk)$ provides

$$|D_k^{-1}b(k)|^2 \leq \mathfrak{s}^2 |D_{k,\text{diag}}^{-1}b(k)|^2 \leq \mathfrak{s}^2 C_{\mathfrak{u}_0} \overline{b}_k^2 v_k^{-1} \leq \mathfrak{s}^2 C_{\mathfrak{u}_0} C_b^2$$

and the assertion follows. □

Combination of the results of Theorem 3.8 and Corollary 3.7 yields the following.

COROLLARY 3.9. *Suppose* (MD) *and* (Dk) *and let the index* $k^*$ *be defined by the balance relation (3.5). Then for* $\Delta = \mathfrak{s}^2 C_{\mathfrak{u}_0} C_b^2$ *and any* $r > 0$

$$\mathbf{E}|v_{k^*}^{-1}(\widehat{\theta} - \widetilde{\theta}_{k^*})^2|^{r/2} \leq \sqrt{e^\Delta \alpha \mathfrak{c}_r} + (2\mathfrak{z}_{k^*})^{r/2},$$

$$v_{k^*}^{-1/2}\mathbf{E}|\widehat{\theta} - \theta| \leq 2\sqrt{e^\Delta} + \sqrt{2\mathfrak{z}_{k^*}}.$$

We conclude this section by a small discussion about relations between of the oracle result and minimax rate of convergence. Most of the theoretical results in the statistical literature are stated about the asymptotic minimax rate of estimation on the functional classes. See for example Lepski (1990, 1992) and Lepski, Mammen and Spokoiny (1997) for pointwise regression estimation or Goldenshluger (1999) and Goldenshluger and Pereversev (2003) for some results in the linear inverse problem. The rate optimal procedures can be obtained using the bias-variance relation (3.5). An immediate corollary of Theorem 3.8 is that the proposed adaptive estimate that selects one out of the family of the spectral cut-off estimates $\widetilde{\theta}_k$ is rate optimal (up to a logarithmic multiplier) for all such set-ups, because it also achieves the accuracy corresponding to the balance relation. A precise formulation of this result lies beyond the focus of this paper.

3.7. *Application to the "sequence space" model.* This section specifies the general results to the sequence space example considered in Section 1.3. In this case, $\widetilde{\theta}_k = y_1 + \cdots + y_{m_k}$, $v_k = \sigma_1^2 + \cdots + \sigma_{m_k}^2$ with $m_1 > m_2 > \cdots > m_K \geq 1$. We additionally assume that $\sigma_i^2$ are monotonously increasing in $i$. The condition (MD) means in this situation that the indices $m_k$ properly decrease to provide an exponential decrease of the sums $v_k$ in $k$. The next result shows that this condition ensures (Dk).

LEMMA 3.2. *For the model (1.3), the condition* (MD) *implies* (Dk) *with the constant* $\mathfrak{s} = (1 - 1/\mathfrak{u}_0)^{-3/2}$.

The proof is given in the Appendix. The estimate $\widetilde{\theta}_k$ has the bias $b_k = \theta_k - \theta = -\sum_{i=m_k+1}^M \mu_i$. The bias-variance relation (3.5) balances the nondecreasing envelope $\overline{b}_k = \max_{l \leq k} |b_l|$ with the variance $v_k^2$ leading to the oracle choice $k^*$. Corollary 3.9 ensures for the adaptive estimate $\widehat{\theta}$ the accuracy of order $v_{k^*}^{-1/2}$ up to the multiplicative factor $\sqrt{2\mathfrak{z}_{k^*}}$.

**4. Simulation.** This section illustrates the performance of the proposed procedure by means of two simulated examples. The first correspond to a severely ill-posed inverse problem with exponentially increasing variances $\sigma_i^2$ and the second to a regularly ill-posed problem with polynomially increasing values $\sigma_i^2$. We focus on two important features of our procedure: "propagation property" and "adaptivity." The "propagation" property means that the selected index only in very few cases is smaller than the oracle one, that means, the "false alarm" situation, when the procedure stops but the modeling bias is still small, is very rare. The "adaptivity" means that the ratio of the risk of the adaptive estimate to the risk of the oracle one is bounded by some fixed constant.

TABLE 1
*Critical values computed under the null hypothesis from 50000 replications, when $K = 20$ and $(\sigma_i = ((n^{2/n})^i)_{i=1,\ldots,n}$ using the sequential procedure*

| $r$ | $\alpha$ | $\mathfrak{z}_1$ | $\mathfrak{z}_2$ | $\mathfrak{z}_3$ | $\mathfrak{z}_4$ | $\mathfrak{z}_5$ | $\mathfrak{z}_6$ | $\mathfrak{z}_7$ | $\mathfrak{z}_8$ | $\mathfrak{z}_9$ | $\mathfrak{z}_{10}$ | $\mathfrak{z}_{11}$ | $\mathfrak{z}_{12}$ | $\mathfrak{z}_{13}$ | $\mathfrak{z}_{14}$ | $\mathfrak{z}_{15}$ | $\mathfrak{z}_{16}$ | $\mathfrak{z}_{17}$ | $\mathfrak{z}_{18}$ | $\mathfrak{z}_{19}$ |
|---|---|---|---|---|---|---|---|---|---|---|---|---|---|---|---|---|---|---|---|---|
| 0.5 | 1.0 | 15.5 | 13.0 | 12.8 | 12.2 | 11.5 | 11.3 | 10.9 | 9.8 | 9.2 | 8.6 | 8.3 | 7.6 | 7.0 | 6.6 | 5.9 | 5.2 | 4.5 | 3.6 | 2.5 |
| 1.0 | 1.0 | 22.5 | 19.0 | 16.4 | 17.2 | 16.2 | 15.6 | 16.8 | 14.4 | 13.4 | 13.2 | 12.9 | 11.9 | 10.2 | 9.3 | 8.3 | 7.3 | 5.8 | 4.7 | 3.4 |

For simplicity we consider "sequence space" models, that is, the data $Y_i$ are generated by the following model: $Y_i = \mu_i + \sigma_i \delta \epsilon_i$, for $i = 1, \ldots, n$ for $n = 50$ and we assume that $\epsilon_i$ are i.i.d. standard normal. In each example the values $(\mu_i)_{i=1,\ldots,n}$ are generated randomly from a centered Gaussian with a decreasing variance $i^{-3}$ and we consider 10 different models of this type. The error level $\delta$ is equal to $10^{-4}, 10^{-5}$ or $10^{-6}$. In every example, the target is the sum of the parameters $\mu_i$, that is, $\theta = \sum_{i=1}^{n} \mu_i$. This set-up is friendly advised by F. Bauer, see for example, Bauer (2007).

We apply the proposed procedure to the family of "weak" estimates $\widetilde{\theta}_k = \sum_{i=1}^{m_k} Y_i$. Our default choice of the "metaparameters" $\alpha$ and $r$ is $\alpha = 1$ and $r = 1/2$. We also report the similar results for $r = 1$, which illustrate that the critical values slightly increase with $r$. More numerical results (not reported here) indicate that the critical values increase with $r$ and decrease with $\alpha$; however, the final results are rather insensitive to the choice of these metaparameters.

In the first example we choose $\sigma_i = a^i$ for $i = 1, \ldots, n$, where $a = n^{2/n}$. We consider the estimates $\widetilde{\theta}_k = \sum_{i=1}^{m_k} Y_i$ with $m_k = [n - 2 * (k - 1)]$, for $k = 1, \ldots, K$ and $K = 20$, then $m_K = 12$.

The critical values $\mathfrak{z}_k$ are computed from 50,000 Monte Carlo replications from the null hypothesis (pure noise model) using the sequential procedure from Section 2.2, see Table 1.

Figure 1 compares the results for our adaptive estimate with the oracle one. The oracle value $k^*$ is defined as $\max\{k : \Delta_k < 1\}$. The results for other values of $\Delta$, for example, $\Delta = 0.5$ or $\Delta = 2$ are very similar and we do not report them here. Each row corresponds to a different level of the noise $\delta$. The panel (a) draws the ratio of the adaptive risk $\mathbf{E}|\widehat{\theta} - \theta|$ obtained from 500 realizations to the corresponding oracle risk $\mathbf{E}|\widetilde{\theta}_{k^*} - \theta|$ for the 10 different

TABLE 2
*Critical values computed under the null hypothesis from* 50000 *replications, when* $K = 15$ *and* $(\sigma_i = i^2)_{i=1,\ldots,n}$ *using the sequential procedure*

| $r$ | $\alpha$ | $\mathfrak{z}_1$ | $\mathfrak{z}_2$ | $\mathfrak{z}_3$ | $\mathfrak{z}_4$ | $\mathfrak{z}_5$ | $\mathfrak{z}_6$ | $\mathfrak{z}_7$ | $\mathfrak{z}_8$ | $\mathfrak{z}_9$ | $\mathfrak{z}_{10}$ | $\mathfrak{z}_{11}$ | $\mathfrak{z}_{12}$ | $\mathfrak{z}_{13}$ | $\mathfrak{z}_{14}$ |
|---|---|---|---|---|---|---|---|---|---|---|---|---|---|---|---|
| 0.5 | 1.0 | 5.5 | 5.0 | 4.6 | 4.3 | 4.1 | 3.9 | 3.4 | 3.1 | 2.8 | 2.6 | 2.2 | 1.7 | 1.3 | 0.9 |
| 1.0 | 1.0 | 8.1 | 7.9 | 6.4 | 6.6 | 7 | 5.8 | 4.8 | 4.3 | 3.9 | 3.6 | 3.0 | 2.0 | 1.5 | 1.0 |

models. In the panel (b) we show the box-plot of $\widehat{k}$ from 500 replications and the "oracles" values $k^*$ (triangles) for the 10 different models described above. One can see that the adaptive risk is in the most of cases not more than twice larger than the the oracle risk. The oracle choice $k^*$ is usually smaller than the adaptively selected $\widehat{k}$, which illustrates "propagation" property: procedure does not stop until $k^*$. It is also worth noticing that both the oracle choice $k^*$ and the adaptive values $\widehat{k}$ decrease with the noise, that is, the smaller the noise, the more coefficients $y_i$ are taken for estimating the sum $\theta = \sum_i \mu_i$.

In the second example we consider a model with $(\sigma_i = i^2)_{i=1,\ldots,n}$ and apply the estimates $\widetilde{\theta}_k = \sum_{i=1}^{m_k} Y_i$ with $m_k = [n/(2^{1/5})^{k-1}]$, for $k = 1, \ldots, K$ and $K = 15$, leading to $m_K = 7$. The critical values $\mathfrak{z}_k$ are computed from 50,000 Monte Carlo replications under the null hypothesis, see Table 2.

Figure 2 presents the results comparing the performance of the adaptive and oracle estimates in the second example. The set-up is the same as in the first example and the results are very similar.

We conclude from this simulation study that the performance of the method is completely in agreement with the theoretical conclusions and the procedure demonstrates quite reasonable performance in all the examples including regular and severely ill-posed problems and for different configurations of the signal and different noise levels.

## APPENDIX.

We start with some useful technical result. Let $(\xi_1, \xi_2)$ be a Gaussian vector with zero mean, $\mathbf{E}\xi_1^2 = \mathbf{E}\xi_2^2 = 1$ and $\rho = \mathbf{E}\xi_1\xi_2$. The correlation coefficient $\rho$ uniquely describes the joint distribution of $\xi_1$ and $\xi_2$ enabling to define for $r \geq 0$ and $\mathfrak{z} \geq 0$

$$Q_r(\rho, \mathfrak{z}) \stackrel{\text{def}}{=} \mathbf{E}[|\xi_1|^{2r}\mathbf{1}(\xi_2^2/2 > \mathfrak{z})], \qquad Q_r^*(\mathfrak{z}) \stackrel{\text{def}}{=} \sup_\rho Q(\rho, \mathfrak{z}).$$

Below we utilize some simple bounds on the quantities $Q_r(\rho, \mathfrak{z})$ and $Q_r^*(\mathfrak{z})$.

LEMMA A.1. *For any $r > 0$ and any $\mathfrak{z} \geq 1$*

$$Q_r^*(\mathfrak{z}) \leq \{C_1(r) + C_2(r)\mathfrak{z}^r\}\mathfrak{z}^{-1/2}e^{-\mathfrak{z}},$$

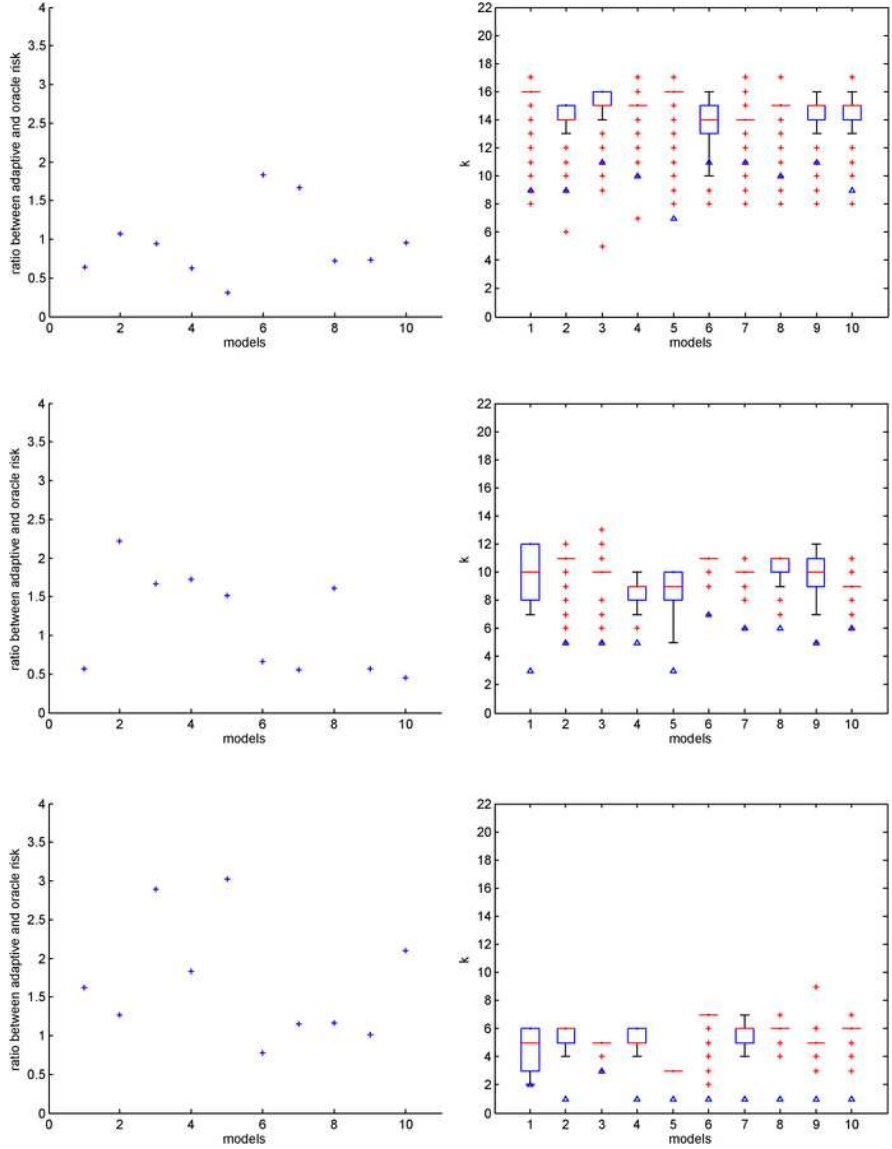

FIG. 1. *The result for the first example with $\delta = 10^{-4}$ (top), $\delta = 10^{-5}$ (middle) and $\delta = 10^{-6}$ (bottom). Left: the ratio of the adaptive risk $\mathbf{E}|\widehat{\theta} - \theta|$ to the oracle risk $\mathbf{E}|\widetilde{\theta}_{k^*} - \theta|$ as function of the model. Right: the boxplots of the adaptive index $\widehat{k}$ based on 500 runs. The triangles show the oracle values $k^*$.*

where $C_1(r)$ and $C_2(r)$ depend on $r$ only. Moreover, for any $\mathfrak{z} \geq 1$

$$\inf_\rho Q_r(\rho, \mathfrak{z}) \geq C_3(r) \mathfrak{z}^{-1/2} e^{-\mathfrak{z}}.$$

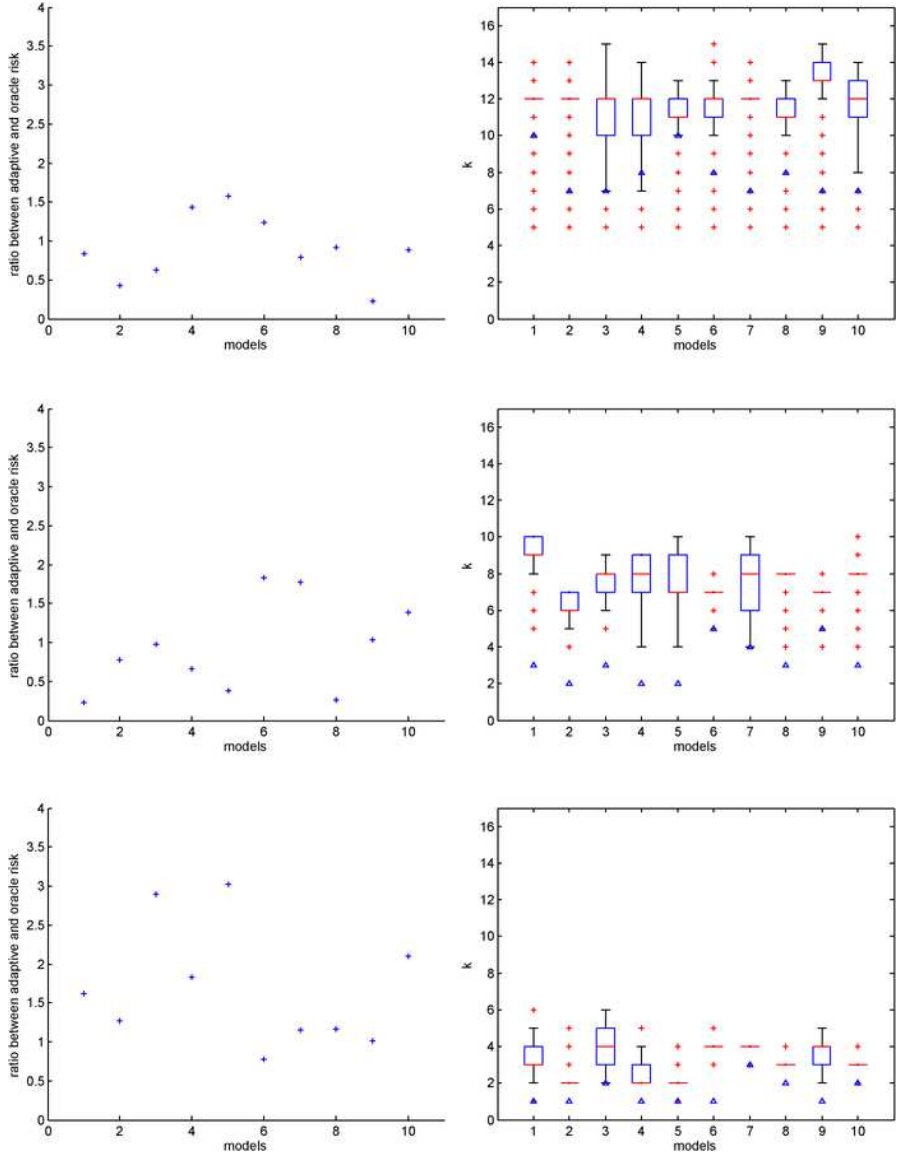

Fig. 2. *The result for the second example with $\delta = 10^{-4}$ (top), $\delta = 10^{-5}$ (middle) and $\delta = 10^{-6}$ (bottom). Left: the ratio of the adaptive risk $\mathbf{E}|\widehat{\theta} - \theta|$ to the oracle risk $\mathbf{E}|\widetilde{\theta}_{k^*} - \theta|$ as function of the model. Right: the boxplots of the adaptive index $\widehat{k}$ based on 500 runs. The triangles show the oracle values $k^*$.*

PROOF. Represent $\xi_1$ as $\rho\xi_2 + \widetilde{\rho}\widetilde{\xi}_1$ where $\widetilde{\rho}$ fulfills $\rho^2 + \widetilde{\rho}^2 = 1$ and $\widetilde{\xi}_1$ is standard normal and independent of $\xi_2$. Note that

$$Q_r(\rho, \mathfrak{z}) = \mathbf{E}|\rho\xi_2 + \widetilde{\rho}\widetilde{\xi}_1|^{2r}\mathbf{1}(\xi_2^2/2 > \mathfrak{z})$$

$$= 0.5\mathbf{E}|\rho\xi_2 + \widetilde{\rho}\widetilde{\xi}_1|^{2r}\mathbf{1}(\xi_2^2/2 > \mathfrak{z}) + 0.5\mathbf{E}|\rho\xi_2 - \widetilde{\rho}\widetilde{\xi}_1|^{2r}\mathbf{1}(\xi_2^2/2 > \mathfrak{z}).$$

One can easily see that there are constants $c_1(r), c_1'(r) > 0$ such that for any $x, y$

$$c_1'(r)\{|\rho x|^{2r} + |\widetilde{\rho} y|^{2r}\} \leq |\rho x + \widetilde{\rho} y|^{2r} + |\rho x - \widetilde{\rho} y|^{2r} \leq c_1(r)\{|\rho x|^{2r} + |\widetilde{\rho} y|^{2r}\}.$$

It is straightforward to check that for some other constants $0 < c_2'(r) < c_2(r)$, $0 < c_3'(r) < c_3(r)$ and $\mathfrak{z} \geq 1$

$$c_2'(r)\mathfrak{z}^{r-1/2}e^{-\mathfrak{z}} \leq \mathbf{E}|\xi_2|^{2r}\mathbf{1}(\xi_2^2 > 2\mathfrak{z}) \leq c_2(r)\mathfrak{z}^{r-1/2}e^{-\mathfrak{z}},$$

$$c_3'(r)\mathfrak{z}^{-1/2}e^{-\mathfrak{z}} \leq \mathbf{E}\mathbf{1}(\xi_2^2 > 2\mathfrak{z}) \leq c_3(r)\mathfrak{z}^{-1/2}e^{-\mathfrak{z}}.$$

The simple algebra yields now

$$Q_r(\rho, \mathfrak{z}) \leq 0.5 c_1(r)\mathbf{E}\{|\rho\xi_2|^{2r} + |\widetilde{\rho}\widetilde{\xi}_1|^{2r}\}\mathbf{1}(\xi_2^2 > 2\mathfrak{z})$$

$$\leq 0.5 c_1(r)\{c_2(r)\mathfrak{z}^r + c_3(r)\mathfrak{c}_r\}\mathfrak{z}^{-1/2}e^{-\mathfrak{z}},$$

$$Q_r(\rho, \mathfrak{z}) \geq 0.5 c_1'(r)\mathbf{E}\{|\rho\xi_2|^{2r} + |\widetilde{\rho}\widetilde{\xi}_1|^{2r}\}\mathbf{1}(\xi_2^2/2 > \mathfrak{z}) \geq C_3(r)\mathfrak{z}^{-1/2}e^{-\mathfrak{z}}.$$

as required. □

PROOF OF THEOREM 3.1. Define for every $m < k \leq K$ the random set $\mathcal{B}_{mk} \stackrel{\text{def}}{=} \{\widehat{\theta}_k = \widetilde{\theta}_m\}$. The definition of the procedure implies

$$\mathcal{B}_{mk} \subseteq \bigcup_{l=1}^{m} \mathbf{1}(v_l^{-1}(\widetilde{\theta}_l - \widetilde{\theta}_k)^2/2 > \mathfrak{z}_l)$$

and

$$\mathbf{E}_0 |v_k^{-1}(\widetilde{\theta}_k - \widehat{\theta}_k)^2|^r \mathbf{1}(\mathcal{B}_{mk}) \leq \sum_{l=1}^{m} \mathbf{E}_0 |v_k^{-1}(\widetilde{\theta}_k - \widetilde{\theta}_m)^2|^r \mathbf{1}(v_l^{-1}(\widetilde{\theta}_l - \widetilde{\theta}_k)^2/2 > \mathfrak{z}_l).$$

Define for $l < m \leq k$

$$v_{lm} = \operatorname{Var}(\widetilde{\theta}_l - \widetilde{\theta}_m), \qquad \xi_{lm} \stackrel{\text{def}}{=} (\widetilde{\theta}_l - \widetilde{\theta}_m)/v_{lm}^{1/2}, \qquad \rho_{lmk} \stackrel{\text{def}}{=} \mathbf{E}_0 \xi_{lk}\xi_{mk}.$$

The conditions of the theorem imply that $v_{lm} \leq \gamma v_l$ for all $l < m$. Therefore

$$\mathbf{E}_0 |v_k^{-1}(\widetilde{\theta}_k - \widehat{\theta}_k)^2|^r \mathbf{1}(\mathcal{B}_{mk}) \leq \sum_{l=1}^{m} \mathbf{E}_0 \left|\frac{\gamma v_m}{v_k}\right|^r |\xi_{mk}|^{2r}\mathbf{1}(\xi_{lk}^2/2 > \mathfrak{z}_l/\gamma)$$

and

$$\mathbf{E}|v_k^{-1}(\widetilde{\theta}_k - \widehat{\theta}_k)^2|^r = \sum_{m=1}^{k-1} \mathbf{E}_0 |v_k^{-1}(\widetilde{\theta}_k - \widehat{\theta}_k)^2|^r \mathbf{1}(\mathcal{B}_{mk})$$

$$\leq \sum_{m=1}^{k-1}\sum_{l=1}^{m} \left|\frac{\gamma v_m}{v_k}\right|^r Q_r(\rho_{lmk}, \mathfrak{z}_l/\gamma)$$

$$\leq \gamma^r \sum_{l=1}^{k-1} Q_r^*(\mathfrak{z}_l/\gamma) \sum_{m=l}^{k-1}\left|\frac{v_m}{v_k}\right|^r.$$

Condition $(MD)$ implies that

$$\sum_{m=l}^{k-1}\left|\frac{v_m}{v_k}\right|^r \le \left|\frac{v_l}{v_k}\right|^r \sum_{m=l}^{k-1} \mathfrak{u}_0^{-(m-l)} \le C(\mathfrak{u}_0)\left|\frac{v_l}{v_k}\right|^r,$$

where $C(\mathfrak{u}_0) = (1 - \mathfrak{u}_0^{-1})^{-1}$. This and Lemma A.1 yield

$$\mathbf{E}|v_k^{-1}(\widetilde{\theta}_k - \widehat{\theta}_k)^2|^r \le \gamma^r C(\mathfrak{u}_0) \sum_{l=1}^{k-1} Q_r^*(\mathfrak{z}_l/\gamma)\left|\frac{v_l}{v_k}\right|^r$$

$$\le C(r, \gamma, \mathfrak{u}_0) \sum_{l=1}^{k-1} \exp\{-\mathfrak{z}_l/\gamma + r\log(v_l/v_k) + r\log(\mathfrak{z}_l)\}.$$

and it remains to check that the choice $\mathfrak{z}_l = a_1 \log(K) + \gamma \log(\alpha^{-1}) + r\gamma \times \log(v_l/v_K)$ with a properly selected $a_1 = a_1(r, \gamma, \mathfrak{u}_0, \mathfrak{u})$ provides in view of $(k-l)\log(\mathfrak{u}_0) \le \log(v_l/v_k) \le (k-l)\log(\mathfrak{u})$ the required bound $\mathbf{E}_0|v_k^{-1}(\widetilde{\theta}_k - \widehat{\theta}_k)^2|^r \le \alpha\mathfrak{c}_r$ for all $k \le K$ and Theorem 3.1 follows. $\square$

PROOF OF THEOREM 3.2. We use again the decomposition

$$\mathbf{E}_0|v_K^{-1}(\widehat{\theta} - \widetilde{\theta}_K)^2|^r = \sum_{k=1}^{K-1} \mathbf{E}_0|v_K^{-1}(\widetilde{\theta}_k - \widetilde{\theta}_K)^2|^r \mathbf{1}(\widehat{k} = k)$$

$$\ge \mathbf{E}_0|v_K^{-1}(\widetilde{\theta}_k - \widetilde{\theta}_K)^2|^r \mathbf{1}(\widehat{k} = k)$$

for any $k < K$. The definition of $\widehat{k}$ implies in the considered case with $\mathfrak{z}_1 = \cdots = \mathfrak{z}_{k-1} = \infty$ that

$$\mathbf{1}(\widehat{k} = k) = \mathbf{1}(v_k^{-1}(\widetilde{\theta}_{k+1} - \widetilde{\theta}_k)^2/2 > \mathfrak{z}_k).$$

With $\rho = \rho_{k,k+1,K} = \mathbf{E}_0 \xi_{k,k+1}\xi_{k,K}$ it holds

$$\mathbf{E}_0|v_K^{-1}(\widehat{\theta} - \widetilde{\theta}_K)^2|^r \ge \mathbf{E}_0|v_K^{-1}(\widetilde{\theta}_k - \widetilde{\theta}_K)^2|^r \mathbf{1}(v_k^{-1}(\widetilde{\theta}_{k+1} - \widetilde{\theta}_k)^2/2 > \mathfrak{z}_k)$$

$$= (v_{k,K}/v_K)^r \mathbf{E}_0|\xi_{k,K}|^{2r} \mathbf{1}(\xi_{k,k+1}^2/2 > \mathfrak{z}_k v_k/v_{k,k+1})$$

$$= (v_{k,K}/v_K)^r Q_r(\rho, \mathfrak{z}_k v_k/v_{k,k+1}).$$

The propagation condition (2.4) implies now that

$$\log(\alpha\mathfrak{c}_r) \ge r\log(v_{k,K}/v_K) + \log Q_r(\rho, \mathfrak{z}_k v_k/v_{k,k+1})$$

yielding in view of Lemma A.1 that

$$\mathfrak{z}_k \ge \frac{v_{k,k+1}}{v_k}\{r\log(v_{k,K}/v_K) + \log\alpha^{-1} - \text{Const.}\log(1 + \log(v_{k,K}/v_K))\}$$

with some fixed constant Const. depending on $r$ only. $\square$

PROOF OF LEMMA 3.2. It suffices to show that the minimal eigenvalue of the matrix $M_k = D_{k,\text{diag}}^{-1} B_k D_{k,\text{diag}}^{-1}$ is bounded away from zero, or, equivalently, the largest eigenvalue of $M_k^{-1}$ is bounded from above: $\|M_k^{-1}\|_\infty \leq (1-1/\mathfrak{u}_0)^{-3}$. Clearly $\mathbf{E}_0 \widetilde{\theta}_j \widetilde{\theta}_l = \mathbf{E}_0 \widetilde{\theta}_l^2 = v_l$ for $j \leq l$, and $M_k$ is the symmetric matrix composed by the elements of the form $\rho_{jl} = v_j^{-1/2} v_l^{-1/2} \mathbf{E}_0 \widetilde{\theta}_j \widetilde{\theta}_l = (v_j/v_l)^{1/2}$ for $j \leq l$. In other words, $M_k$ is the covariance matrix for the set of random variables $\eta_l = \widetilde{\theta}_l / v_l^{1/2}$ for $l = 1, \ldots, k$.

Define $\gamma_l = v_l^{-1/2}(\widetilde{\theta}_l - \widetilde{\theta}_{l+1})$ for $l < k$ and $\gamma_k = v_k^{-1/2} \eta_k$. The random variables $\gamma_l$ are independent zero mean normal with the variance $s_l \stackrel{\text{def}}{=} \mathbf{E} \gamma_l^2 = v_l^{-1}(v_l - v_{l+1})$ for $l < k$ and $s_k = 1$. The condition $(MD)$ implies for all $l \leq k$ that $(1 - 1/\mathfrak{u}_0) \leq s_l \leq (1 - 1/\mathfrak{u})$. Define $\gamma^{(k)} = (\gamma_1, \ldots, \gamma_k)^\top$ and $\eta^{(k)} = (\eta_1, \ldots, \eta_k)^\top$. The identities $\gamma_l = \eta_l - \eta_{l+1}(v_{l+1}/v_l)^{1/2}$ for $l < k$ can be written as $\gamma^{(k)} = A_k \eta^{(k)}$, where line $l$ of the matrix $A_k$ only contains only two nonzero entries: $a_{l,l} = 1$ and $a_{l,l+1} = -v_{l+1}^{1/2}/v_l^{1/2}$ for $l = 1, \ldots, k-1$. Again, the condition $(MD)$ implies that $\|I - A_k\|_\infty \leq 1/\mathfrak{u}_0$ and $\|A_k^{-1}\|_\infty = \|\{I - (I - A_k)\}^{-1}\| \leq (1 - 1/\mathfrak{u}_0)^{-1}$. Similar bound holds for $A_k^\top$. Obviously $\mathbf{E}_0 \gamma^{(k)}(\gamma^{(k)})^\top = \varGamma_k \stackrel{\text{def}}{=} \text{diag}(s_1, \ldots, s_k)$. This yields

$$\varGamma_k = \mathbf{E} A_k \eta^{(k)}(\eta^{(k)})^\top A_k^\top = A_k M_k A_k^\top$$

and $\|M_k^{-1}\|_\infty \leq \|A_k^{-1}\|_\infty^2 \cdot \|\varGamma_k^{-1}\|_\infty \leq (1 - 1/\mathfrak{u}_0)^{-3}$, then the result follows. □

**Acknowledgments** The authors thank the associate editor and two anonymous referees for very useful remarks and suggestions leading to a serious improvement of the paper.

WEIERSTRASS-INSTITUTE AND
HUMBOLDT UNIVERSITY BERLIN
MOHRENSTR. 39
10117 BERLIN
GERMANY
E-MAIL: spokoiny@wias-berlin.de

MODAL'X AND UNIVERSITÉ
EUROPÉENNE DE BRETAGNE
CAMPUS DE KER LANN
RUE BLAISE PASCAL-BP 37203
35172 BRUZ CEDEX
FRANCE
E-MAIL: celine.vial@univ-rennes1.fr